\pdfoutput=1
\RequirePackage{ifpdf}
\ifpdf 
\documentclass[pdftex]{sigma}
\else
\documentclass{sigma}
\fi

\usepackage{tikz-cd}

\newcommand{\Asf}{\mathsf{A}}
\newcommand{\Csf}{\mathsf{C}}

\newcommand{\Xfk}{\mathfrak{X}}

\newcommand{\Frm}{\mathrm{F}}
\newcommand{\irm}{\mathrm{i}}
\newcommand{\prm}{\mathrm{p}}

\newcommand{\Nat}{\mathbb{N}}
\newcommand{\R}{\mathbb{R}}

\newcommand{\util}{\tilde{u}}
\newcommand{\xtil}{\tilde{x}}

\newcommand{\abs}[1]{\lvert #1 \rvert}
\renewcommand{\dh}{{\rm d}_\mathrm{h}}
\newcommand{\dhn}{{\rm d}_{\mathrm{h}\n}}
\newcommand{\dv}{{\rm d}_\mathrm{v}}
\newcommand{\n}{\nabla}
\newcommand{\p}{\partial}
\newcommand{\Ptil}{\tilde{P}}
\newcommand{\Phat}{\hat{P}}
\newcommand{\Stil}{\tilde{S}}
\newcommand{\Shat}{\hat{S}}

\newcommand{\pd}[2]{\frac{\p #1}{\p #2}}
\newcommand{\vf}[1]{\frac{\p}{\p #1}}
\newcommand{\vfe}[2]{\vf{#1}\biggr|_{#2}}

\newcommand{\pdm}[3]{\frac{\p^{\abs{#3}} #1}{\p #2^{#3}}}

\newcommand{\hook}{\text{$\,\mathop{\vrule height0.6pt width3pt depth0pt\kern0pt\vrule height6pt width0.6pt depth0pt}\,\,$}}

\newcommand{\Jpi}{J^1\pi}
\newcommand{\Jbpi}{J^2\pi}
\newcommand{\Jbpihat}{\widehat{J}^2\pi}
\newcommand{\Jbkpi}{J^{2k}\pi}
\newcommand{\Jbkapi}{J^{2k-1}\pi}
\newcommand{\Jkpi}{J^k\pi}
\newcommand{\Jkpihat}{\widehat{J}^k\pi}
\newcommand{\Jkapi}{J^{k-1}\pi}
\newcommand{\JkApi}{J^{k+1}\pi}
\newcommand{\Jlpi}{J^l\pi}
\newcommand{\Jftypi}{J^\infty \pi}
\newcommand{\Jpia}{J^1\pia}
\newcommand{\Jpika}{J^1\pika}

\newcommand{\Jotau}{J^0\tau}
\newcommand{\Jtau}{J^1\tau}
\newcommand{\Jtaua}{J^1\taua}
\newcommand{\Jbtau}{J^2\tau}
\newcommand{\Jbtauhat}{\widehat{J}^2\tau}
\newcommand{\Jltau}{J^l\tau}
\newcommand{\Jlatau}{J^{l-1}\tau}

\newcommand{\Jltauhat}{\widehat{J}^l\tau}
\newcommand{\Jtaul}{J^1\taul}
\newcommand{\Jtaula}{J^1\taula}

\newcommand{\JlAtauhat}{\widehat{J}^{l+1}\tau}

\newcommand{\jpphi}{j^1_p\phi}
\newcommand{\jbpphi}{j^2_p\phi}
\newcommand{\jkpphi}{j^k_p\phi}
\newcommand{\jkapphi}{j^{k-1}_p\phi}

\newcommand{\jkppsi}{j^k_p\psi}

\newcommand{\jpoG}{j^1_p\omega}

\newcommand{\jpzG}{j^1_p\zeta}
\newcommand{\jbpzG}{j^2_p\zeta}
\newcommand{\jkapzG}{j^{k-1}_p\zeta}

\newcommand{\jGG}{j^1\Gamma}

\newcommand{\jtaulla}{j^1\taulla}

\newcommand{\pia}{\pi_1}
\newcommand{\piao}{\pi_{1,0}}
\newcommand{\piaa}{(\pia)_1}
\newcommand{\piaao}{(\pia)_{1,0}}
\newcommand{\pib}{\pi_2}
\newcommand{\pibo}{\pi_{2,0}}
\newcommand{\piba}{\pi_{2,1}}
\newcommand{\pik}{\pi_k}
\newcommand{\piko}{\pi_{k,0}}
\newcommand{\pikl}{\pi_{k,l}}
\newcommand{\pilk}{\pi_{l,k}}
\newcommand{\pika}{\pi_{k-1}}
\newcommand{\pikka}{\pi_{k,k-1}}
\newcommand{\pikaa}{(\pika)_1}
\newcommand{\pikaao}{(\pika)_{1,0}}
\newcommand{\pifty}{\pi_\infty}
\newcommand{\piftyk}{\pi_{\infty,k}}

\newcommand{\taua}{\tau_1}
\newcommand{\tauao}{\tau_{1,0}}
\newcommand{\taul}{\tau_l}
\newcommand{\taulo}{\tau_{l,0}}
\newcommand{\taula}{\tau_{l-1}}

\newcommand{\taulaao}{(\taula)_{1,0}}
\newcommand{\taulao}{\tau_{l-1,0}}
\newcommand{\taulla}{\tau_{l,l-1}}
\newcommand{\tauLao}{(\taul)_{1,0}}

\newcommand{\iaa}{\irm_{1,1}}
\newcommand{\iaka}{\irm_{1,k-1}}

\DeclareMathOperator{\id}{id}

\numberwithin{equation}{section}

\newtheorem{Theorem}{Theorem}[section]
\newtheorem*{Theorem*}{Theorem}

\newtheorem{Proposition}[Theorem]{Proposition}
 { \theoremstyle{definition}

 }

\begin{document}

\allowdisplaybreaks

\renewcommand{\thefootnote}{}

\newcommand{\arXivNumber}{2309.01594}

\renewcommand{\PaperNumber}{013}

\FirstPageHeading

\ShortArticleName{Lepage Equivalents and the Variational Bicomplex}

\ArticleName{Lepage Equivalents and the Variational Bicomplex\footnote{This paper is a~contribution to the Special Issue on Symmetry, Invariants, and their Applications in honor of Peter J.~Olver. The~full collection is available at \href{https://www.emis.de/journals/SIGMA/Olver.html}{https://www.emis.de/journals/SIGMA/Olver.html}}}

\Author{David SAUNDERS}

\AuthorNameForHeading{D.~Saunders}

\Address{Lepage Research Institute, Pre\v{s}ov, Slovakia}
\Email{\href{mailto:david@symplectic.email}{david@symplectic.email}}

\ArticleDates{Received October 05, 2023, in final form January 30, 2024; Published online February 09, 2024}

\Abstract{We show how to construct, for a Lagrangian of arbitrary order, a Lepage equivalent satisfying the closure property:\ that the Lepage equivalent vanishes precisely when the Lagrangian is null. The construction uses a homotopy operator for the horizontal differential of the variational bicomplex. A choice of symmetric linear connection on the manifold of independent variables, and a global homotopy operator constructed using that connection, may then be used to extend any global Lepage equivalent to one satisfying the closure property. In the second part of the paper we investigate the r\^{o}le of vertical endomorphisms in constructing such Lepage equivalents. These endomorphisms may be used directly to construct local homotopy operators. Together with a symmetric linear connection they may also be used to construct global vertical tensors, and these define infinitesimal nonholonomic projections which in turn may be used to construct Lepage equivalents. We conjecture that these global vertical tensors may also be used to define global homotopy operators.}


\Classification{58A10; 58A20; 83D05}

\renewcommand{\thefootnote}{\arabic{footnote}}
\setcounter{footnote}{0}

\subsection*{Dedication}

In the notes to Section~5 of~\cite{Olv86}, Peter Olver wrote about the variational complex and the variational bicomplex `It is hoped that these methods will inspire further research in the geometric theory of the calculus of variations'. A few years later~\cite{Olv93} he wrote `In the geometric theory of the calculus of variations in mechanics, the Cartan form, which first arose as the integrand in Hilbert's invariant integral, plays a ubiquitous role'. Lepage equivalents are generalizations of Cartan forms, and I hope that this paper will be a small contribution to Peter's project.

\section{Introduction}

In recent years there has been a revival of interest in the `fundamental Lepage equivalent' of a~Lagrangian, a~differential form on a jet bundle which (as with any such Lepage equivalent) provides a geometrical construction leading to the Euler--Lagrange equations of the corresponding variational problem, but which has the additional property that it is closed precisely when the Lagrangian is null~\cite{Pal22,Urb22}. The original formulation of the fundamental Lepage equivalent was given for first order Lagrangians (in~\cite{Kru77}, and then independently in~\cite{Bet84}). Although expressed in local coordinates, the form is in fact invariant under changes of coordinates, and so is a global geometric object. There had, however, been no similar construction for higher order Lagrangians.

A construction for Lagrangians of arbitrary order has now been proposed in~\cite{Voi22}, giving a~Lepage form of order no greater than $4k-2$ for a Lagrangian of order $k$. The construction is again given in local coordinates, but now there is no guarantee that it will be defined globally. In addition, if the original Lagrangian happens to be first order, the new Lepage form will in general be of second order and will differ from the original, first order, fundamental Lepage equivalent.

In the first part of this paper, after giving some background on the different types of Lepage equivalent, we propose a new method of constructing a fundamental Lepage equivalent for a~Lagrangian of arbitrary order by using homotopy operators for the horizontal differential in the variational bicomplex on the infinite jet manifold. This has the disadvantage that any bound on the order of the resulting Lepage form, although necessarily finite, will depend on the number of independent variables. On the other hand, the choice of a symmetric linear connection will allow the construction of a global form, and in the case of a first order Lagrangian the result will be independent of any connection and we recover the classic fundamental Lepage equivalent.

In the second part of the paper, we explore the potential for clarifying this construction by using `vertical endomorphisms' on jet bundles, tensorial objects depending on a closed differential form, which are related to the canonical isomorphism between the tangent space at any point of an affine space, and the vector space on which the affine space is modelled~\cite{Sau87}. We recall how local homotopy operators for the horizontal differential can be constructed from these vertical endomorphisms, and then we show how a symmetric linear connection can be used to remove the dependence on the differential form to produce a globally-defined, fully tensorial object. (A related but technically different approach has been described in~\cite{Bet93}.) Such a `vertical tensor' can be used to give an infinitesimal rigidity to nonholonomic jet bundles, allowing the construction of a global Lepage equivalent for a Lagrangian of arbitrary order. Finally, we offer a conjecture regarding how these vertical tensors, together with a covariant version of the horizontal differential, might be used to construct a global homotopy operator for the ordinary horizontal differential.

\section{Notation}

We adopt a modified version of the notation used in~\cite{Sau89}. We let $\pi \colon E \to M$ be a fibred manifold with $\dim M = m$ and $\dim E = m + n$. The $k$-jet manifold of $\pi$ will be denoted $\Jkpi$ with projections $\pik \colon \Jkpi \to M$, $\piko \colon \Jkpi \to E$ and $\pikl \colon \Jkpi \to \Jlpi$ where $l < k$. A typical element of $\Jkpi$ will be denoted $\jkpphi$. We use similar notation for jets of the cotangent bundle $\tau \colon T^* M \to M$.
We let $\Xfk\bigl(\Jkpi\bigr)$ denote the module of vector fields on $\Jkpi$, and $\Omega^r\bigl(\Jkpi\bigr)$ the module of $r$-forms.

Regarding the jet bundle $\pika \colon \Jkapi \to M$ as the starting bundle, we shall let $\pikaa$: $\Jpika \to M$ denote its first jet bundle, and we let $\iaka \colon \Jkpi \to \Jpika$ be the canonical inclusion. There is also an intermediate submanifold $\Jkpihat \subset \Jpika$, the semiholonomic manifold~\cite{Lib97}, with a~canonical symmetrization projection $\prm_k \colon \Jkpihat \to \Jkpi$.

We also use the infinite jet bundle $\pifty \colon \Jftypi \to M$ where $\Jftypi$, although infinite dimensional, is a Fr\'{e}chet manifold and so is reasonably well behaved. We let $\Omega^r$ (without specifying a~manifold) denote the module of $r$-forms on $\Jftypi$ of globally finite order, so that each such form is projectable to some $\Jkpi$.

Any differential form $\omega \in \Omega^r$ can be decomposed uniquely into its contact components
\[
\omega = \omega^{(0)} + \omega^{(1)} + \cdots + \omega^{(p)} + \cdots + \omega^{(r)},
\]
where if $r > m$ then $\omega^{(p)} = 0$ for $p < r-m$. We let $\Omega^{p,q} \subset \Omega^r$ with $p+q=r$ denote the submodule of $p$-contact $r$-forms. In a similar way, a differential form $\omega \in \Omega^r\bigl(\Jkpi\bigr)$ on a finite order jet manifold may be decomposed into contact components, but these will normally be defined on $\JkApi$ rather than on $\Jkpi$.

When using coordinates, we take fibred coordinates $\bigl(x^i,u^\alpha\bigr)$ on $E$ over base coordinates $\bigl(x^i\bigr)$ on $M$. Jet coordinates will be denoted $(u^\alpha_i)$ on $\Jpi$ and $\bigl(u^\alpha_i, u^\alpha_{(ij)}\bigr)$ on $\Jbpi$ with parentheses denoting symmetrization because $u^\alpha_{(ji)}$ is the same coordinate as $u^\alpha_{(ij)}$. For this reason we use the symbol $\#(ij)$ to equal $1$ when $i=j$ and to equal $2$ when $i \ne j$, in order to avoid double counting during summation.

On higher order jet manifolds this notation becomes unwieldy and we write $(u^\alpha_I)$ instead, where $I \in \Nat^m$ is a multi-index with $I(i)$ giving the number of copies of the index $i$, so that this notation automatically takes care of symmetrization. We let $1_i$ denote the multi-index with a~single $1$ in the $i$-th position; $\abs{I} = \sum_{i=1}^m I(i)$ is the length of $I$, and $I! = I(1)! I(2)! \cdots I(m)!$ is its factorial. Any summation involving multi-indices will be indicated explicitly, including the zero multi-index where appropriate.

Sometimes we need to use a mixed notation, and converting to or from multi-index notation requires coefficients to be adjusted. If $F(J)$ is some object depending on the multi-index $J$, then
\[
\sum_{\abs{J} = r+1} \frac{\abs{J}!}{J!} F(J) = \sum_{i=1}^m \sum_{\abs{I} = r} \frac{\abs{I}!}{I!} F(I+1_i),
\]
where the quotient $\abs{J}! / J!$ is the `weight' of the multi-index $J$.

We use notation
\[
\theta^\alpha = {\rm d}u^\alpha - u^\alpha_j {\rm d}x^j , \qquad \theta^\alpha_i = {\rm d}u^\alpha_i - u^\alpha_{(ij)} {\rm d}x^j , \qquad \theta^\alpha_I = {\rm d}u^\alpha_I - u^\alpha_{I+1_j} {\rm d}x^j
\]
for local contact $1$-forms, and
\[
\omega_0 = {\rm d}x^1 \wedge \cdots \wedge {\rm d}x^m , \qquad \omega_i = i_{\p / \p x^i} \, \omega_0 = (-1)^{i-1} {\rm d}x^1 \wedge \cdots \wedge \widehat{{\rm d}x^i} \wedge \cdots \wedge {\rm d}x^m
\]
(where the circumflex indicates an omitted factor) for local forms horizontal over $M$. Local total derivatives, dual to the local contact forms, will be denoted $d_i$ and are given explicitly as
\[
\vf{x^i} + u^\alpha_i \vf{u^\alpha} , \qquad \vf{x^i} + u^\alpha_i \vf{u^\alpha_i} + u^\alpha_{(ij)} \vf{u^\alpha_j} , \qquad
\vf{x^i} + \sum_I u^\alpha_{I+1_i} \vf{u^\alpha_I} .
\]
In the finite order case they are vector fields along the map $\pikka$ rather than on a single jet manifold. We also use the symbol $\p_i$ to indicate $\p / \p x^i$ as a vector field along the map $\pik$.

On a nonholonomic jet manifold we need to distinguish between the two levels of jet coordinates, and we use juxtaposition, with a dot to indicate when a particular index is missing. So on $\Jpia$ the coordinates are $\bigl(x^i, u^\alpha_{\cdot\cdot}, u^\alpha_{i \cdot}, u^\alpha_{\cdot j}, u^\alpha_{ij}\bigr)$ and on $\Jpika$ they are $\bigl(x^i, u^\alpha_{I\cdot}, u^\alpha_{Ij}\bigr)$.

Finally, we note that $\piao \colon \Jpi \to E$ is an affine bundle, modelled on the vector bundle $\pi^* T^* M \otimes V\pi$, so that the vertical bundle $V\piao$ is canonically isomorphic to $\pia^* T^* M \otimes \piao^* V\pi$; the inverse of this isomorphism may be regarded as a tensor field
\[
S = \p_i \otimes \theta^\alpha \otimes \vf{u^\alpha_i} ,
\]
a section of the bundle $ \pia^* TM \otimes T^* \Jpi \otimes T\Jpi$ over $\Jpi$. We shall call this the \emph{first order vertical tensor}.

\section{Background}

Many of the results in the geometrical calculus of variations can be described in terms of \emph{source forms} and \emph{Lepage equivalents} (see~\cite{KKS10} for a useful summary and historical references).

A source form is a form $\varepsilon \in \Omega^{m+1}\bigl(\Jlpi\bigr)$ with the properties that it is horizontal over $E$, and maximally horizontal over $M$, so that in coordinates it appears as $\varepsilon = \varepsilon_\alpha \theta^\alpha \wedge \omega_0$. The zero set of a source form is a submanifold of $\Jlpi$ representing a family of partial differential equations; if $\lambda = L \omega_0 \in \Omega^m\bigl(\Jkpi\bigr)$ is a horizontal $m$-form, a Lagrangian, then it gives rise to a source form $\varepsilon_\lambda \in \Omega^{m+1}\bigl(\Jbkpi\bigr)$, the \emph{Euler--Lagrange form} of $\lambda$, incorporating the Euler--Lagrange equations of the variational problem defined by $\lambda$:
\[
\varepsilon_\lambda = \sum_{\abs{I}=0}^k (-1)^{\abs{I}} \, d_I \biggl( \pd{L}{u^\alpha_I} \biggr) \theta^\alpha \wedge \omega_0 .
\]
A Lepage form is a form $\vartheta \in \Omega^m\bigl(\Jlpi\bigr)$ with the property that $({\rm d}\vartheta)^{(1)}$, the $1$-contact component of its exterior derivative, is a source form. A Lepage equivalent of a Lagrangian $\lambda \in \Omega^m\bigl(\Jkpi\bigr)$ is a Lepage form $\vartheta_\lambda \in \Omega^m\bigl(\Jlpi\bigr)$ with $l \ge k$ such that the difference $\vartheta_\lambda - \pilk^* \lambda$ is a contact form; the corresponding source form $({\rm d}\vartheta_\lambda)^{(1)}$ is then just the Euler--Lagrange form $\varepsilon_\lambda$. Different Lepage equivalents of the same Lagrangian give the same Euler--Lagrange form.

If $m = 1$, then each Lagrangian $\lambda$ gives rise to a unique globally-defined Lepage equivalent $\vartheta_\lambda \in \Omega^1\bigl(\Jbkapi\bigr)$, the Cartan form of the Lagrangian. However, complications arise when $m \ge 2$, and these concern both existence and uniqueness. Clearly if $\vartheta_\lambda$ is a Lepage equivalent of~$\lambda$, then so is $\vartheta_\lambda + {\rm d}\psi + \omega$ where $\omega$ is at least $2$-contact, and in fact the converse is true: if~$\vartheta_\lambda$,~$\vartheta_\lambda^\prime$ are both Lepage equivalents of $\lambda$, then locally $\theta_\lambda^\prime - \theta_\lambda = {\rm d}\psi + \omega$, where $\psi$ will necessarily be a~contact form and may be chosen to be $1$-contact.

As far as existence is concerned, if we initially consider forms which are at most $1$-contact, then locally
\begin{equation}
\label{EpLe}
\vartheta_\lambda = L \omega_0 + \sum_{\abs{J}=0}^{k-1} \sum_{\abs{K}=0}^{k-\abs{J}-1}
\frac{(-1)^{\abs{J}} (J+K+1_j)! \abs{J}! \abs{K}!}{(\abs{J}+\abs{K}+1)! J! K!}\, d_J
\biggl( \pd{L}{u^\alpha_{J+K+1_j}} \biggr) \theta^\alpha_K \wedge \omega_j
\end{equation}
is a Lepage equivalent, known as the \emph{principal Lepage equivalent} of $\lambda$. When $k=1$ this is just the Poincar\'{e}--Cartan form of $\lambda$,
\[
\vartheta_\lambda = L \omega_0 + \pd{L}{u^\alpha_j} \theta^\alpha \wedge \omega_j
\]
and is defined globally; it is the unique Lepage equivalent of $\lambda$ which is both at most $1$-contact and also horizontal over $E$. When $k=2$, we obtain
\[
\vartheta_\lambda = L \omega_0 + \biggl( \biggl( \pd{L}{u^\alpha} - \frac{1}{\#(ij)} \, d_i \biggl( \pd{L}{u^\alpha_{(ij)}} \biggr) \biggr) \theta^\alpha
+ \frac{1}{\#(ij)} \pd{L}{u^\alpha_{(ij)}} \theta^\alpha_i \biggr) \wedge \omega_j,
\]
which again, perhaps surprisingly, is invariant under a fibred change of coordinates $\xtil = \xtil^j\bigl(x^i\bigr)$, $\util^\beta = \util^\beta\bigl(x^i, u^\alpha\bigr)$ and is therefore also defined globally. For $k \ge 3$, however, there is no such Lepage equivalent invariant under coordinate changes~\cite{HK83}, and so choices need to be made in order to obtain a globally defined form. Several authors (see, for instance, \cite{FF83,GM83}) have used connections of various kinds for this purpose; an approach with an algebraic flavour which uses a~symmetric linear connection on $M$ has been described in~\cite{AB99}, and when the connection is defined locally as the canonical connection arising from a system of coordinates on $M$, then the resulting local Poincar\'{e}--Cartan form is the principal Lepage equivalent~\eqref{EpLe}.

A different approach~\cite{Sau87} has been to `pretend' that the $k$-th order Lagrangian is really first order by using a tubular neighbourhood of $\Jkpi$ in $\Jpika$ and `spreading out' the Lagrangian using the neighbourhood's projection. By repeating this process, a global Lepage equivalent may be constructed. In fact only infinitesimal projections are needed, mapping $T_{\Jkpi}\Jpika$ to $T\Jkpi$, and we shall see in Section~\ref{Sproj} that a symmetric linear connection on $M$ determines a suitable family of projections. In the second order case, it may be seen that only the restriction of the projection to the semiholonomic manifold $\Jbpihat$ is needed, explaining why the symmetrization projection $\prm_2$ may be used to give a global Lepage equivalent in this case.

The Lepage equivalents described so far have all been at most $1$-contact. There have, however, been important examples of Lepage equivalents involving higher contact terms. One such, defined for a nonvanishing first order Lagrangian, is the \emph{Carath\'{e}odory form}~\cite{Car29}
\[
\vartheta_\lambda = \frac{1}{L^{m-1}} \bigwedge_{j=1}^m \biggl( L \, {\rm d}x^j + \pd{L}{u^\alpha_j} \theta^\alpha \biggr) ;
\]
this decomposable form is again defined globally and indeed is invariant, not just under a fibred change of coordinates, but under a general change $\xtil = \xtil^j\bigl(x^i, u^\alpha\bigr)$, $\util^\beta = \util^\beta\bigl(x^i, u^\alpha\bigr)$. A similar form for a nonvanishing second order Lagrangian,
\[
\vartheta_\lambda = \frac{1}{L^{m-1}} \bigwedge_{j=1}^m \biggl( L \, {\rm d}x^j + \biggl( \pd{L}{u^\alpha_j}
- \frac{1}{\#(ij)}\, d_i \biggl( \pd{L}{u^\alpha_{(ij)}} \biggr) \biggr) \theta^\alpha
+ \frac{1}{\#(ij)} \pd{L}{u^\alpha_{(ij)}} \theta^\alpha_i \biggr),
\]
was described in~\cite{Olv93} (see also~\cite{CS04}).

The Lepage equivalent of particular interest in the present paper, again involving higher contact terms, is the \emph{fundamental Lepage equivalent} of a first order Lagrangian~\cite{Bet84, Kru77}
\[
\vartheta_\lambda = \sum_{p=0}^{\min\{m,n\}} \frac{1}{(p!)^2} \frac{\p^p L}{\p u^{\alpha_1}_{j_1} \cdots \p u^{\alpha_p}_{j_p}} \theta^{\alpha_1} \wedge \cdots \wedge \theta^{\alpha_p} \wedge \omega_{j_1 \dots j_p} .
\]
This satisfies the \emph{closure property}, that ${\rm d}\vartheta_\lambda = 0$ precisely when the Lagrangian is null:\ that is, when the Euler--Lagrange form $\varepsilon_\lambda$ is zero. (Of course any individual form $\vartheta_\lambda$ is either closed or not closed; the closure property applies to the procedure mapping $\lambda$ to $\vartheta_\lambda.$) Once again this form is defined globally, and indeed is invariant under a general change of coordinates $\xtil = \xtil^j\bigl(x^i, u^\alpha\bigr)$, $\util^\beta = \util^\beta\bigl(x^i, u^\alpha\bigr)$~\cite{CS05}. The content of the closure property lies in the requirement that ${\rm d}\vartheta_\lambda = 0$ when $\lambda$ is null; for any Lepage equivalent $\vartheta_\lambda$ it is obvious that the converse holds, that $\lambda$ is null when ${\rm d}\vartheta_\lambda = 0$.
In the next section, we consider how the construction of the fundamental Lepage equivalent might be generalised for higher order Lagrangians.

\section{The closure property}

The question of whether it is possible to find a procedure for constructing a Lepage equivalent which satisfies the closure property, although solved in 1977 for first order Lagrangians, has been an open problem for higher order Lagrangians (see~\cite{Pal22,Urb22} and the references therein). An original solution to this problem was given in~\cite{Voi22}, using the Vainberg--Tonti Lagrangian of a~source form $\varepsilon$, the horizontal $m$-form $\lambda_\varepsilon$ obtained locally in coordinates from $\varepsilon = \varepsilon_\alpha \theta^\alpha \wedge \omega_0$ by the fibred homotopy operator
\[
\lambda_\varepsilon = u^\alpha \int_0^1 \varepsilon_\alpha\bigl(x^i, tu^\beta_I\bigr) {\rm d}t .
\]
Typically $\lambda_\varepsilon$ has the same order as $\varepsilon$. If in fact $\varepsilon = \varepsilon_\lambda$, so that the source form is the Euler--Lagrange form of a given Lagrangian $\lambda$, then the Vainberg--Tonti Lagrangian $\lambda_{\varepsilon_\lambda}$ and the pullback of $\lambda$ have the same Euler--Lagrange equations so that they differ by $h({\rm d}\alpha)$ for some horizontal $(m-1)$-form $\alpha$. Then, taking $\vartheta_{\lambda_{\varepsilon_\lambda}}$ to be the principal Lepage equivalent of the Vainberg--Tonti Lagrangian in the given coordinates and writing $\vartheta^\Frm = \vartheta_{\lambda_{\varepsilon_\lambda}}+ {\rm d}\alpha$ we find that, to within pullbacks,
\[
\bigl({\rm d}\vartheta^\Frm\bigr)^{(1)} = {\rm d}\vartheta_{\lambda_{\varepsilon_\lambda}}^{(1)} = \varepsilon_{\lambda_{\varepsilon_\lambda}} = \varepsilon_\lambda
\]
so that ${\rm d}\vartheta^\Frm$ is a source form, and that
\[
h\bigl(\vartheta^\Frm\bigr) = h\bigl(\vartheta_{\lambda_{\varepsilon_\lambda}}\bigr) + h({\rm d}\alpha) = \lambda_{\varepsilon_\lambda} + h({\rm d}\alpha) = \lambda
\]
so that $\vartheta^\Frm$ is a Lepage equivalent of $\lambda$, and finally that
\[
{\rm d}\vartheta^\Frm = {\rm d}\vartheta_{\lambda_{\varepsilon_\lambda}},
\]
so that if $\lambda$ is a null Lagrangian then $\varepsilon_\lambda = 0$ and therefore ${\rm d}\vartheta^\Frm = 0$.

This procedure therefore satisfies the closure property. It is not, though, a generalisation of the fundamental Lepage equivalent for first order Lagrangians, because it is always at most $1$-contact, whereas the fundamental Lepage equivalent is obtained from the Poincar\'{e}--Cartan form by adding higher contact terms.

We shall, instead, define an alternative procedure which uses homotopy operators for the horizontal differential of the variational bicomplex to add the higher contact terms. Recall that this bicomplex is defined for forms of globally finite order on the infinite jet manifold $\Jftypi$, as shown in the diagram, see Figure~\ref{Fig1} (note that the squares with vertical arrows labelled $\pifty^*$ commute, whereas those with vertical arrows labelled $\dv$ anticommute). The rows and columns are all locally exact, and indeed all the $\dh$ rows apart from the first are globally exact~\cite{And89, Tak79, Tuj82, Tul80, Vin78}. Any Lagrangian $\lambda \in \Omega^m\bigl(\Jkpi\bigr)$ will have a pullback $\piftyk^* \lambda \in \Omega^{0,m}$ on $\Jftypi$ which for simplicity we shall continue to denote by $\lambda$ without the pullback map.
\begin{figure}[ht]\centering
$
\begin{tikzcd}
& & 0 \arrow{d} & 0 \arrow{d} & 0 \arrow{d} && 0 \arrow{d} \\
0 \arrow{r} & \R \arrow{r} & \Omega^0(M) \arrow{r}{{\rm d}} \arrow{d}{\pifty^*} & \Omega^1(M) \arrow{r}{{\rm d}} \arrow{d}{\pifty^*}
& \Omega^2(M) \arrow{d}{\pifty^*} \arrow{r}{{\rm d}} & \cdots \arrow{r}{{\rm d}} & \Omega^m(M) \arrow{d}{\pifty^*} \\
0 \arrow{r} & \R \arrow{r} & \Omega^{0,0} \arrow{r}{\dh} \arrow{d}{\dv} & \Omega^{0,1} \arrow{r}{\dh} \arrow{d}{\dv}
& \Omega^{0,2} \arrow{d}{\dv} \arrow{r}{\dh} & \cdots \arrow{r}{\dh} & \Omega^{0,m} \arrow{d}{\dv} \\
& 0 \arrow{r} & \Omega^{1,0} \arrow{r}{\dh} \arrow{d}{\dv} & \Omega^{1,1} \arrow{r}{\dh} \arrow{d}{\dv}
& \Omega^{1,2} \arrow{d}{\dv} \arrow{r}{\dh} & \cdots \arrow{r}{\dh} & \Omega^{1,m} \arrow{d}{\dv} \\
& 0 \arrow{r} & \Omega^{2,0} \arrow{r}{\dh} \arrow{d}{\dv} & \Omega^{2,1} \arrow{r}{\dh} \arrow{d}{\dv}
& \Omega^{2,2} \arrow{d}{\dv}
\arrow{r}{\dh} & \cdots \arrow{r}{\dh} & \Omega^{2,m} \arrow{d}{\dv} \\
& & \vdots & \vdots & \vdots && \vdots
\end{tikzcd}
$
\caption{The variational bicomplex.}\label{Fig1}
\end{figure}

Let $\vartheta_\lambda$ denote the pullback to $\Jftypi$ of any local Lepage equivalent of $\lambda$ which is at most $1$-contact, so that $\vartheta_\lambda^{(1)} = \vartheta_\lambda - \lambda \in \Omega^{1,m-1}$, and let $P$ denote any local homotopy operator for the $\dh$ rows (apart from the first) of the variational bicomplex. Define the \emph{extension of $\vartheta_\lambda$ by $P$} to be the $m$-form defined locally by
\begin{align*}
\begin{split}
\vartheta^\Frm & = \vartheta_\lambda + \sum_{p=1}^{m-1} (-P\dv)^p \vartheta_\lambda^{(1)}
 = \lambda + \vartheta_\lambda^{(1)} - (P\dv) \vartheta_\lambda^{(1)} + (P\dv)^2 \vartheta_\lambda^{(1)} - \cdots + (-P\dv)^{m-1} \vartheta_\lambda^{(1)} \\
& \in \Omega^{0,m} \oplus \Omega^{1,m-1} \oplus \Omega^{2,m-2} \oplus \Omega^{3,m-3} \oplus \cdots \oplus \Omega^{m,0} ,
\end{split}
\end{align*}
so that $\vartheta^\Frm$ is another Lepage equivalent of $\lambda$. We shall show that this method of constructing~$\vartheta^\Frm$ satisfies the closure property, by diagram chasing.

Suppose that $\lambda$ is a null Lagrangian, so that
\[
0 = \varepsilon_\lambda = \bigl({\rm d}\vartheta^\Frm\bigr)^{(1)} = \dv \vartheta_\lambda^{(0)} + \dh \vartheta_\lambda^{(1)} = \dv \lambda + \dh \vartheta_\lambda^{(1)} .
\]
Then
\begin{align*}
\bigl({\rm d}\vartheta^\Frm\bigr)^{(2)} & = \dv \bigl(\vartheta^{\Frm(1)}\bigr) + \dh \bigl(\vartheta^{\Frm(2)}\bigr)
 = \dv \vartheta_\lambda^{(1)} - \dh P \dv \vartheta_\lambda^{(1)} \\
& = P \dh \dv \vartheta_\lambda^{(1)} = - P \dv \dh \vartheta_\lambda^{(1)} = P \dv\dv\lambda = 0
\end{align*}
using the homotopy property $\dh \circ P + P \circ \dh = \id$, and in a similar way
\begin{align*}
\bigl({\rm d}\vartheta^\Frm\bigr)^{(p+1)} & = \dv \bigl(\vartheta^{\Frm(p)}\bigr) + \dh \bigl(\vartheta^{\Frm(p+1)}\bigr)
 = \dv \bigl(\vartheta^{\Frm(p)}\bigr) - \dh P \dv \bigl(\vartheta^{\Frm(p)}\bigr) \\
& = P \dh \dv \bigl(\vartheta^{\Frm(p)}\bigr) = - P \dv \dh \bigl(\vartheta^{\Frm(p)}\bigr) = P \dv \dv \bigl(\vartheta^{\Frm(p-1)}\bigr) = 0
\end{align*} for $2 \le p \le m-1$, where the penultimate equality arises recursively from
\[
\dh \bigl(\vartheta^{\Frm(p)}\bigr) + \dv \bigl(\vartheta^{\Frm(p-1)}\bigr) = \bigl({\rm d}\vartheta^{\Frm}\bigr)^{(p)} = 0 .
\]
Thus ${\rm d}\vartheta^\Frm = \bigl({\rm d}\vartheta^\Frm\bigr)^{(m+1)}$, and we may see that this maximal contact component also vanishes by traversing the diagram in the opposite direction. For $2 \le p \le m$, we have
\begin{align*}
\dh \dv \bigl(\vartheta^{\Frm(p)}\bigr) & = - \dh \dv P \dv \bigl(\vartheta^{\Frm(p-1)}\bigr) = \dv \dh P \dv \bigl(\vartheta^{\Frm(p-1)}\bigr) \\
& = \dv \dv \bigl(\vartheta^{\Frm(p-1)}\bigr) - (\dv P) \dh \dv \bigl(\vartheta^{\Frm(p-1)}\bigr)
= - (\dv P) \dh \dv \bigl(\vartheta^{\Frm(p-1)}\bigr) ,
\end{align*}
but
\[
\dh \dv \bigl(\vartheta^{\Frm(1)}\bigr) = \dh \dv \vartheta_\lambda^{(1)} = - \dv \dh \vartheta_\lambda^{(1)} = \dv \dv \lambda = 0
\]
so that $\dh \dv (\vartheta^{\Frm(m)}) = 0$. As $\dh \colon \Omega^{m+1,0} \to \Omega^{m+1,1}$ is injective by exactness, we see finally that $\bigl({\rm d}\vartheta^\Frm\bigr)^{(m+1)} = \dv \bigl(\vartheta^{\Frm(m)}\bigr) = 0$. We shall describe suitable local homotopy operators for $\dh$, constructed using vertical endomorphisms, in the next Section; by using them we obtain the following result.

\begin{Theorem}
Let $\lambda$ be the pullback to $\Jftypi$ of a Lagrangian of any order, and let $\vartheta_\lambda$ be the pullback to $\Jftypi$ of any local Lepage equivalent of $\lambda$ which is at most $1$-contact. A local homotopy operator $P$ then defines a local Lepage equivalent $\vartheta^\Frm$ which is an extension of $\vartheta_\lambda$ and which satisfies the closure property, that ${\rm d}\vartheta^\Frm = 0$ precisely when $\lambda$ is null.
\end{Theorem}

We can also consider a global version of this result, noting that the diagram chasing above would apply equally to global operators as it does to local ones. We have seen that additional structures, such as connections or nonholonomic projections, are needed to construct a global Lepage equivalent when the order of the Lagrangian is greater than two. A global homotopy operator for the horizontal differential on $\Jftypi$ has also been found~\cite[Theorem 5.56]{And89} and again this uses a symmetric linear connection on the base manifold $M$.

\begin{Theorem}
Let $\lambda$ be the pullback to $\Jftypi$ of a Lagrangian of any order, and let $\vartheta_\lambda$ be the pullback to $\Jftypi$ of any global Lepage equivalent of $\lambda$ which is at most $1$-contact, constructed using additional data as appropriate. A global homotopy operator $P$, such as the one described in~{\rm \cite{And89}} using a symmetric linear connection, then defines a global Lepage equivalent $\vartheta^\Frm$ which is an extension of $\vartheta_\lambda$ and which satisfies the closure property, that ${\rm d}\vartheta^\Frm = 0$ precisely when $\lambda$ is null.
\end{Theorem}

We remark that in fact there is no requirement for the homotopy operators in each term to be the same, and we could generalise the formula to
\[
\vartheta^\Frm = \vartheta_\lambda + \sum_{p=1}^{m-1} (-1)^p (P_{p+1}\dv P_p \dv \cdots P_2 \dv) \vartheta_\lambda^{(1)},
\]
where $P_p$ is a homotopy operator for the $p$-contact row of the variational bicomplex.


\section{Vertical endomorphisms}

The most basic example of a `vertical endomorphism' is the almost tangent structure on a tangent manifold $TM$. This is simply a tensorial expression of the isomorphism between a vector space and its tangent space at any point, applied to the tangent spaces to a manifold, and may be regarded as a $1$-form taking values in the sub-bundle of $TTM \to TM$ containing the vertical vectors. A similar object may be defined using a more complicated procedure on a higher order tangent manifold $T^k M$~\cite{CCS86}, now giving a $1$-form taking its values in the sub-bundle of vertical vectors in $TT^k M \to T^k M$.

Vertical endomorphisms $S^\eta$ on jet manifolds $\Jkpi$, where $\eta \in \Omega^1(M)$ is a closed $1$-form, were defined in~\cite{Sau87}. The construction started with a point $\jkpphi \in \Jkpi$ and a tangent vector $\xi$ at $\jkapphi \in \Jkapi$ vertical over $M$. Any such vector may be represented by a $1$-parameter family of local sections $\phi_t$ where $\phi_0 = \phi$ and $\xi$ is the tangent vector at $t=0$ to the curve $t \mapsto \jkapphi_t$. Given a function $f$ on $M$ defined in a neighbourhood of $p$, the \emph{vertical lift} of $\xi$ to $\jkpphi$ in the direction specified by ${\rm d}f$ then used the $1$-parameter family of local sections $\psi_t \colon q \mapsto \phi_{tf(q)}(q)$ to define a~curve $\jkppsi_t$ in $\Jkpi$ and therefore a tangent vector at $\jkpphi$. The vertical endomorphism~$S^\eta$ at any point $\jkpphi \in \Jkpi$ was then defined by starting with any tangent vector in $T_{\jkpphi}\Jkpi$, projecting it to $T_{\jkapphi}\Jkapi$, taking the vertical representative using the contact structure, and then applying the vertical lift (using any function $f$ satisfying $f(p)=0$ and ${\rm d}f = \eta$ in a neighbourhood of $p$) to give a new tangent vector in $T_{\jkpphi}\Jkpi$. It may be shown that this construction is well defined, and so independent of the choices of $\phi_t$ and $f$, and that it gives a tensor field $S^\eta \in \Omega^1\bigl(\Jkpi\bigr) \otimes \Xfk\bigl(\Jkpi\bigr)$ expressed in coordinates\footnote{The numerical coefficient given in~\cite[equation~(3.4)]{Sau87} and repeated in~\cite[after Definition~6.5.6]{Sau89} is incorrect as it does not take account of the use of an individual index~$i$ in a multi-index formula.} as
\[
S^\eta = \sum_{\abs{J}+\abs{K} \le k-1} \frac{(J+K+1_i)!}{J!K!(\abs{K}+1)} \pdm{\eta_i}{x}{K} \, \theta^\alpha_J \otimes \vf{u^\alpha_{J+K+1_i}} .
\]
It is evident from this formula that, when acting on forms, the operators $S^\eta$ on $\Jkpi$ and on $\Jlpi$ with $l>k$ are related by the pullback map $\pilk^*$, so that we may define a similar operator acting on forms on $\Jftypi$ without ambiguity.

Given local coordinates $\bigl(x^i\bigr)$ on $U \subset M$, we write $S^i$ rather than $S^{{\rm d}x^i}$ for the operators on $U^\infty = \pifty^{-1}(U)$. These local operators have the rather simpler coordinate description
\[
S^i = \sum_{\abs{I}=0}^\infty \bigl( I(i)+1 \bigr) \theta^\alpha_I \otimes \vf{u^\alpha_{I+1_i}}
\]
and may be used to construct local homotopy operators for the horizontal differential~$\dh$ on~$\Jftypi$. One such homotopy operator, involving an ordering of the coordinates $x^i$, was given in~\cite{Tul80}. Other homotopy operators, not using such an ordering, may be constructed from two different repeated actions of $S^i$ on forms: these are
\[
\Stil^J = i_{S^{j_1} \circ S^{j_2} \circ \cdots \circ S^{j_r}} , \qquad \Shat^J = i_{S^{j_1}} \circ i_{S^{j_2}} \circ \cdots \circ i_{S^{j_r}},
\]
where $\abs{J}=r$ and $J = 1_{j_1} + 1_{j_2} + \cdots + 1_{j_r}$. Note that the first action is a derivation, whereas the second is not if $r>1$; the multi-index notation is justified because operators $S^i$ and $S^j$ commute. The following result was obtained in~\cite[Theorem~1]{CS09}.

\begin{Proposition}
Define the differential operators $\Ptil, \Phat \colon \Omega^{p,q}(U^\infty) \to \Omega^{p,q-1}(U^\infty)$, with $p \ge1$ and $1 \le q \le m$, by $\Ptil(\omega) = i_{{\rm d}/{\rm d}x^i} \bigl( \Ptil^i(\omega) \bigr)$, $\Phat(\omega) = i_{{\rm d}/{\rm d}x^i} \bigl( \Phat^i(\omega) \bigr)$, where
\begin{align*}
\Ptil^i(\omega)
& = \sum_{I=0}^\infty \frac{(-1)^{\abs{I}}(m-q)! \abs{I}!}{p(m-q+\abs{I}+1)! I!} \, d_I \Stil^{I+1_i} \omega , \\
\Phat^i(\omega)
& = \sum_{I=0}^\infty \frac{(-1)^{\abs{I}}(m-q)! \abs{I}!}{p^{\abs{I}+1}(m-q+\abs{I}+1)! I!} \, d_I \Shat^{I+1_i} \omega .
\end{align*}
Then both $\Ptil$ and $\Phat$ are homotopy operators for $\dh$.
\end{Proposition}

In general the operators $\Ptil$ and $\Phat$ are different, although they are equal when acting on forms projectable to $\Jpi$, and also when acting on forms in $\Omega^{1,q}(U^\infty)$. In the latter case, writing the operator as $P$, \cite[Theorem~2]{CS09} gives
\[
\omega - \dh P \omega = \theta^\alpha \wedge \sum_{\abs{I}=0}^\infty (-1)^{\abs{I}} \, d_I \bigl( i_{\p / \p u^\alpha_I} \omega \bigr)
\]
for any $\omega \in \Omega^{1,m}$, so that $\omega - \dh P \omega$ is a source form. If we write $\vartheta_\lambda = \lambda - P \dv \lambda$, then
\[
({\rm d}\vartheta_\lambda)^{(1)} = \dv \bigl(\vartheta_\lambda^{(0)}\bigr) + \dh \bigl(\vartheta_\lambda^{(1)}\bigr) = (\dv \lambda) - \dh P (\dv \lambda)
\]
so that in particular $({\rm d}\vartheta_\lambda)^{(1)}$ is a source form. Thus $\vartheta_\lambda$ is a Lepage form, and it is clearly a local Lepage equivalent of $\lambda$. In coordinates with $\lambda = L \omega_0$
\[
\vartheta_\lambda = L \omega_0 + \sum_{\abs{J}, \abs{K} = 0}^\infty
\frac{(-1)^{\abs{J}}(J+K+1_j)! \abs{J}! \abs{K}!}{(\abs{J}+\abs{K}+1)! J! K!}
\, d_J \biggl( \pd{L}{u^\alpha_{J+K+1_j}} \biggr) \theta_K \wedge \omega_j ,
\]
so that it is the pullback to $\Jftypi$ of the principal Lepage equivalent of $\lambda$. We obtain the above formula from those for $S$ and $P$ by using the multi-index Leibniz' rule and the identity for weighted sums of binomial coefficients
\[
\sum_{0 \le K \le I} \frac{(-1)^{\abs{K}}I!}{(\abs{K}+p+1) K! (I-K)!} = \frac{p! \abs{I}!}{(\abs{I}+p+1)!}
\]
obtained by first evaluating the integral $\int_0^1 x^p (x-1)^r {\rm d}x$ in two different ways, and then using the Vandermonde identity for the convolution of scalar binomial coefficients.

We can also use the homotopy operators $P$ to give a simple proof of the result mentioned earlier, that the difference between two Lepage equivalents for the same Lagrangian is the sum of a closed form and a form which is at least $2$-contact. Let $\vartheta_\lambda$ and $\vartheta^\prime_\lambda$ be two Lepage equivalents for the Lagrangian $\lambda$, and put $\vartheta = \vartheta_\lambda - \vartheta^\prime_\lambda$, so that $\vartheta^{(0)} = 0$. As $\vartheta_\lambda$ and $\vartheta^\prime_{\lambda}$ give rise to the same Euler--Lagrange form, we see that
\[
\dh\bigl(\vartheta^{(1)}\bigr) + \dv\bigl(\vartheta^{(0)}\bigr) = ({\rm d}\vartheta)^{(1)} = 0,
\]
so that $\dh\bigl(\vartheta^{(1)}\bigr) = 0$ and therefore locally $\vartheta^{(1)} = \dh P \bigl(\vartheta^{(1)}\bigr)$. Then
\[
\vartheta^{(1)} = {\rm d} P (\vartheta^{(1)}) - \dv P (\vartheta^{(1)}),
\]
where $\dv P \bigl(\vartheta^{(1)}\bigr) \in \Omega^{2,m-2}$, so that
\begin{align*}
\vartheta = {\rm d} P \bigl(\vartheta^{(1)}\bigr) + \bigl( \vartheta^{(2)} - \dv P \bigl(\vartheta^{(1)}\bigr) \bigr) + \cdots + \vartheta^{(m)}
 \in {\rm d}\Omega^{1,m-2} \oplus \Omega^{2,m-2} \oplus \cdots \oplus \Omega^{m,0} .
\end{align*}

Finally, in this Section we apply these operators to a first order Lagrangian $\lambda = L\omega_0$. As $\dv \lambda = {\rm d}\lambda$ is also first order, we see that
\[
\lambda - P\dv \lambda = L \omega_0 + S^i \biggl( \pd{L}{u^\alpha} \theta^\alpha + \pd{L}{u^\alpha_j} \theta^\alpha_j \biggr) \omega_i
= L \omega_0 + \pd{L}{u^\alpha_i} \theta^\alpha \wedge \omega_i ,
\]
the local expression of the Poincar\'{e}--Cartan form. This is also first order, and we then see that each successive operator $P^i$ is simply a multiple of $S^i$. We obtain
\[
(-P\dv)^p \lambda = \frac{1}{(p!)^2} \frac{\p^p L}{\p u^{\alpha_1}_{i_1} \cdots \p u^{\alpha_p}_{i_p}} \theta^{\alpha_1} \wedge \cdots \wedge \theta^{\alpha_p} \wedge \omega_{i_1 \dots i_p} ,
\]
showing that $\sum_{p=0}^m (-P\dv)^p \lambda$ gives the local expression of the standard fundamental Lepage equivalent of a first order Lagrangian.

\section{Connections and vertical tensors}

We have seen that, for a first order Lagrangian, the fundamental Lepage equivalent may be constructed locally using the homotopy operators $S^i$, and also that it is a global object which may be constructed using the first order vertical tensor~$S$. These are two different facets of the same construction and they arise because, in the first order case, the formulation of a vertical endomorphism $S^\eta$ does not in fact require the $1$-form $\eta$ to be closed. From the coordinate description
\[
S^\eta = \eta_i \theta^\alpha \otimes \vf{u^\alpha_i},
\]
it is clear that at any point $\jpphi \in \Jpi$ the value of $S^\eta$ depends only on the cotangent vector $\eta|_p$ and not on the values of $\eta$ at any other points. Thus, given any cotangent vector $\eta|_p \in T^*_p M$, we may choose a closed $1$-form $\zeta$ in a neighbourhood of $p$ satisfying $\zeta|_p = \eta|_p$, for example, the form given in coordinates centred on $p$ by $\zeta = {\rm d}\bigl(\eta_i(p) x^i\bigr)$, and put
\[
S^\eta|_{\jpphi} = S^\zeta|_{\jpphi} \in \bigl(T^* \Jpi \otimes T\Jpi\bigr)_{\jpphi} .
\]
Doing this at each point of $\Jpi$ gives a well-defined vertical endomorphism $S^\eta$ for an arbitrary $1$-form $\eta \in \Omega^1(M)$, and it is clear that the mapping $\eta \mapsto S^\eta$ is just that given by the vertical tensor $S$.

The same approach will not work directly for higher order vertical endomorphisms. For example, the coordinate description of $S^\eta$ on $\Jbpi$ is
\[
S^\eta = \eta_i \theta^\alpha \otimes \vf{u^\alpha_i} + \frac{1}{\#(ij)} \pd{\eta_i}{x^j} \theta^\alpha \otimes \vf{u^\alpha_{(ij)}}
+ \frac{2}{\#(ij)}\eta_i \theta^\alpha_j \otimes \vf{u^\alpha_{(ij)}}
\]
and at any point $\jbpphi \in \Jbpi$ the value of $S^\eta$ depends, not just on the cotangent vector $\eta_p$, but also on the derivative of the $1$-form $\eta$ at $p$.

We can, however, deal with this problem by supposing that we are given a symmetric linear connection $\n$ on $M$; the infinitesimal parallel translation defined by $\n$ will then provide enough information to specify the derivative of $\eta$. The vertical endomorphism defined by the $1$-form $\eta$ (not necessarily closed) and the connection $\n$ will be given in coordinates as
\[
S^\eta_\n = \eta_i \theta^\alpha \otimes \vf{u^\alpha_i} + \frac{1}{\#(hj)} \eta_i \Gamma^i_{hj} \theta^\alpha \otimes \vf{u^\alpha_{(hj)}}
+ \frac{2}{\#(ij)} \eta_i \theta^\alpha_j \otimes \vf{u^\alpha_{(ij)}},
\]
where $\Gamma^i_{hj}$ are the connection coefficients of $\n$, so that the mapping $\eta \to S^\eta_\n$ will define a second order vertical tensor $S_\n$, a section of the bundle $\pib^* TM \otimes T^* \Jbpi \otimes T \Jbpi$ over $\Jbpi$.

Formally, as the $1$-form $\eta$ is a section of the cotangent bundle $\tau \colon T^* M \to M$, we regard the connection $\n$ as a linear Ehresmann connection $\Gamma \colon T^* M \to \Jtau$, a section of the jet bundle $\tauao \colon \Jtau \to T^* M$, so that the connection coefficients $\Gamma^i_{hj}$ are just the jet coordinates of $\Gamma$. \big(Of~course, the connection $\n$ also defines a linear Ehresmann connection on the tangent bundle, but there the jet coordinates are~$-\Gamma^i_{hj}$.\big) For each $\jbpphi \in \Jbpi$, we may choose a~closed $1$-form~$\zeta$ in a~neighbourhood of $p$ satisfying $\jpzG = \Gamma(\eta|_p)$, for example the form given in coordinates centred on $p$ by $\zeta = {\rm d} \bigl( \eta_i(p) x^i + \tfrac{1}{2} \eta_i(p) \Gamma^i_{hj}(p) x^h x^j \bigr)$, and put
\[
S_\n^\eta\big|_{\jbpphi} = S^\zeta\big|_{\jbpphi} \in \bigl(T^* \Jbpi \otimes T\Jbpi\bigr)_{\jbpphi} .
\]
Doing this at each point of $\Jbpi$ now gives a well-defined vertical endomorphism $S_\n^\eta$ for an arbitrary $1$-form $\eta \in \Omega^1(M)$, and so we can construct a second order vertical tensor with coordinate expression
\begin{equation}
\label{Snabla}
S_\n = \p_i \otimes \biggl( \theta^\alpha \otimes \vf{u^\alpha_i} + \frac{1}{\#(hj)} \Gamma^i_{hj} \theta^\alpha \otimes \vf{u^\alpha_{(hj)}}
+ \frac{2}{\#(ij)} \theta^\alpha_j \otimes \vf{u^\alpha_{(ij)}} \biggr) .
\end{equation}

A similar procedure may be carried out for higher order vertical endomorphisms, but requires the use of semiholonomic jets to allow for symmetrization. For example, in the third order case we use the connection map $\Gamma \colon T^* M \to \Jtau$, regarded as a bundle morphism $\tau \to \taua$ over the identity on $M$, and its prolongation $\jGG \colon \Jtau \to \Jtaua$. The composition $\jGG \circ \Gamma$ then takes its values in the semiholonomic manifold $\Jbtauhat \subset \Jtaua$~\cite[Section~5.3]{Sau89}, so that if $\prm_2 \colon \Jbtauhat \to \Jbtau$ is the symmetrization projection, then we may use
\[
\Gamma_2 = \prm_2 \circ \jGG \circ \Gamma \colon\ T^* M \to \Jbtau
\]
as the map which allows us to specify the first and second derivatives at $p$ of the closed local $1$-form $\zeta$ by setting $\jbpzG = \Gamma_2(\eta|_p)$.

More generally, we construct the maps $\Gamma_l$ recursively. Suppose we have the map $\Gamma_{l-1} \colon T^* M \to \Jlatau$, and that it is a section of $\taulao$ with the property that $\jGG_{l-1} \circ \Gamma$ takes its values in the semiholonomic manifold $\Jltauhat \subset \Jtaula$, so that we may set
\[
\Gamma_l = \prm_l \circ \jGG_{l-1} \circ \Gamma \colon\ T^* M \to \Jltau .
\]
We note first that
\begin{align*}
\taulla \circ \Gamma_l & = \taulaao \circ i_{1,l-1} \circ \Gamma_l
 = \taulaao \circ i_{1,l-1} \circ \prm_l \circ \jGG_{l-1} \circ \Gamma \\
& = \taulaao \circ \jGG_{l-1} \circ \Gamma
 = \Gamma_{l-1} \circ \tauao \circ \Gamma
 = \Gamma_{l-1} ,
\end{align*}
so that
\begin{align*}
\taulo \circ \Gamma_l = \taulao \circ \taulla \circ \Gamma_l = \taulao \circ \Gamma_{l-1} = \id_{T^* M},
\end{align*}
and therefore that $\Gamma_l$ is a section of $\taulo \colon \Jltau \to T^* M$. We must also check that $\jGG_l \circ \Gamma$ takes its values in the semiholonomic manifold $\JlAtauhat$, the submanifold of $\Jtaul$ given by equality of the two maps $\jtaulla$ and $i_{1,l-1} \circ \tauLao$ to $\Jtaula$~\cite[Section 6.2]{Sau89}; but at any point $\jpoG \in \Jtau$ we know that
\begin{align*}
\bigl( \jtaulla \circ \jGG_l \bigr) \bigl(\jpoG\bigr) = j^1(\taulla \circ \Gamma_l) \bigl(\jpoG\bigr)
 = j^1\Gamma_{l-1}\bigl(\jpoG\bigr)
\end{align*}
and
\begin{align*}
\bigl( i_{1,l-1} \circ \tauLao \circ \jGG_l \bigr) \bigl(\jpoG\bigr)
 & = \bigl( i_{1,l-1} \circ \tauLao \bigr) \bigl( j^1_p(\Gamma_l \circ \omega) \bigr) \\
 & = i_{1,l-1} ( \Gamma_l (\omega(p)) )
 = \jGG_{l-1} ( \Gamma(\omega(p)) ) ,
\end{align*}
so that if $\jpoG$ is in the image of $\Gamma$, then $\jpoG = \Gamma(\omega(p))$ and
\[
\bigl( \jtaulla \circ \jGG_l \bigr) ( \Gamma(\omega(p)) ) = \jGG_{l-1} ( \Gamma(\omega(p)) )
= \bigl( i_{1,l-1} \circ \tauLao \circ \jGG_l \bigr) ( \Gamma(\omega(p)) )
\]
as required. We may therefore define $\Gamma_{l+1} = \prm_{l+1} \circ \jGG_l \circ \Gamma$ and continue the process. The recursion starts with $l=2$ and $\Gamma_1 = \Gamma \colon T^* M \to \Jtau$, or even degenerately with $l=1$ and $\Gamma_0 = \id_{T^* M} \colon T^* M \to \Jotau = T^* M$.

To find a coordinate expression for these maps, let $\bigl(x^i, y_j\bigr)$ be the coordinates on $T^* M$, so that the jet coordinates on $\Jtau$ are $y_{ij}$ and on $\Jltau$ are $y_{iJ}$. As the connection is linear and symmetric, we see that $y_{ij} \circ \Gamma = y_h \Gamma^h_{ij}$ with $y_{ji} \circ \Gamma = y_{ij} \circ \Gamma$, and in general if $y_{iJ} \circ \Gamma_l = y_h \Gamma^h_{J+1_i}$, then
\[
y_{iJj} \circ \jGG_l = y_h \pd{\Gamma^h_{J+1_i}}{x^j} + y_{hj} \Gamma^h_{J+1_i}
\]
so that
\[
y_{iJj} \circ \jGG_l \circ \Gamma = y_g \biggl( \pd{\Gamma^g_{J+1_i}}{x^j} + \Gamma^g_{hj} \Gamma^h_{J+1_i} \biggr);
\]
the coordinates $y_{iJ+1_j} \circ \Gamma_{l+1}$ may then be obtained by symmetrization. In the degenerate case, the coordinates of $\Gamma_0$ are of course $\Gamma^h_i = \delta^h_i$.

We have, therefore, obtained the following result.

\begin{Theorem}
Let $\pi \colon E \to M$ be a fibred manifold, and let $\n$ be a symmetric linear connection on $M$. On any jet manifold $\Jkpi$ there is a canonical vertical tensor $S_\n$ defined in the following way. If $\eta \in \Omega^1(M)$ and $\jkpphi \in \Jkpi$, let $\zeta$ be any local closed $1$-form defined in a neighbourhood of~$p$ satisfying $\jkapzG = \Gamma_{k-1}(\eta|_p)$ $($for example, a $1$-form defined using a polynomial in coordinates~$x^i$ centred on $p)$ and put $S^\eta_\n\big|_{\jkpphi} = S^\zeta\big|_{\jkpphi}$. Then $S^\eta_\n\big|_{\jkpphi}$ is independent of the choice of $\zeta$. The resulting map $\jkpphi \mapsto S^\eta_\n\big|_{\jkpphi}$ is a vertical endomorphism depending at each point $\jkpphi$ only on the cotangent vector $\eta|_p$ and so defines a vertical tensor $\eta \mapsto S^\eta_\n$.
\end{Theorem}

The coordinate expression of $S_\n$ is
\[
S_\n = \p_h \otimes \sum_{\abs{J}+\abs{K} \le k-1} \frac{(J+K+1_i)!}{J!K!(\abs{K}+1)}
\Gamma^h_{K+1_i} \theta^\alpha_J \otimes \vf{u^\alpha_{J+K+1_i}} ,
\]
and combining the sums over the index $i$ and the multi-index $K$ in the usual way then gives
\[
S_\n = \p_h \otimes \sum_{\substack{\abs{J}+\abs{K} \le k\\ \abs{K}>0}} \frac{(J+K)!}{J!K!}
\Gamma^h_K \theta^\alpha_J \otimes \vf{u^\alpha_{J+K}} .
\]
A similar formula, without an upper bound on the length of the multi-indices, may be used on~$\Jftypi$ for the map $\pifty^* T^* M \otimes T^* \Jftypi \to T\Jftypi$.


\section{Infinitesimal nonholonomic projections}
\label{Sproj}

As mentioned earlier, two possible approaches to defining global Lepage equivalents for higher order Lagrangians involve using either connections, or tubular neighbourhoods of holonomic jet manifolds inside nonholonomic jet manifolds. We have remarked that the latter approach really involves only the infinitesimal projection defined by the tubular neighbourhood at points of the holonomic submanifold, and we can now see that the existence of vertical tensors allows the two approaches to be related:\ a symmetric linear connection on the base manifold will define an infinitesimal nonholonomic projection $T_{\Jkpi}\Jpika \to T\Jkpi$ for $k \ge 2$.

The simplest example is in the second order case, where $\iaa \colon \Jbpi \to \Jpia$ is the canonical inclusion. We start with a point $\jbpphi \in \Jbpi$ and a tangent vector $\xi \in T_{\jbpphi} \Jpia$ which is vertical over $\Jpi$, so that $\xi \in V_{\jbpphi}\piaao$. We then apply the isomorphism
\[
\Asf \colon\ V\piaao \to \piaa^* T^* M \otimes \piaao^* V\pia
\]
arising from the affine structure of $\piaao \colon \Jpia \to \Jpi$ \big(restricted to points of $\Jbpi$\big) and follow this by $S_\n$, giving a map
\[
p_\n = S_\n \circ \Asf \colon\ V_{J^2\pi}(\pia)_{1,0} \to V\pibo,
\]
so that $p_\n(\xi) \in V_{\jbpphi}\pibo$.

There are, of course, many possible extensions of $p_\n$ to a map $T_{\Jbpi}\Jpia \to T\Jpi$; but there is precisely one such extension satisfying the requirement that $p_\n \circ T\iaa = \id_{T\Jbpi}$. We may see this by looking at coordinate representations. At any point $\iaa\bigl(\jbpphi\bigr) \in \Jpia$
\[
\Asf \biggl( \vf{u^\alpha_{\cdot j}} \biggr) = {\rm d}x^j \otimes \vf{u^\alpha} , \qquad
\Asf \biggl( \vf{u^\alpha_{ij}} \biggr) = {\rm d}x^j \otimes \vf{u^\alpha_i} ,
\]
and composing with $S_\n$ as presented in formula~\eqref{Snabla} in the previous section gives
\[
p_\n \biggl( \vf{u^\alpha_{\cdot j}} \biggr) = \vf{u^\alpha_j} + \frac{1}{\#(ik)} \Gamma^j_{ik} \vf{u^\alpha_{(ik)}} , \qquad
p_\n \biggl( \vf{u^\alpha_{ij}} \biggr) = \frac{1}{\#(ij)} \vf{u^\alpha_{(ij)}} .
\]
(Nominally the image of, say,
\[
\vfe{u^\alpha_{ij}}{\iaa(\jbpphi)} \mapsto {\rm d}x^j\big|_p \otimes \vfe{u^\alpha_i}{\jpphi} \in T^*_p M \otimes T_{\jpphi}\Jpi
\]
is not directly in the domain of $S_\n$; but as $S_\n$ incorporates the projection $T\piba \colon T\Jbpi \to T\Jpi$ we may represent that image by an element of $T^*_p M \otimes T_{\jbpphi}\Jbpi$ without ambiguity.) Noting now that
\begin{alignat*}{3}
&T\iaa \biggl( \vf{x^i} \biggr) = \vf{x^i}, \qquad &&
T\iaa \biggl( \vf{u^\alpha_i} \biggr) = \vf{u^\alpha_{i\cdot}} + \vf{u^\alpha_{\cdot i}}, & \\
&T\iaa \biggl( \vf{u^\alpha} \biggr) = \vf{u^\alpha_{\cdot\cdot}}, \qquad &&
T\iaa \biggl( \vf{u^\alpha_{(ij)}} \biggr) = \vf{u^\alpha_{ij}}, &
\end{alignat*}
we see that necessarily we must have
\[
p_\n \biggl( \vf{x^i} \biggr) = \vf{x^i} , \qquad
p_\n \biggl( \vf{u^\alpha_{\cdot\cdot}} \biggr) = \vf{u^\alpha} , \qquad
p_\n \biggl( \vf{u^\alpha_{i\cdot}} \biggr) = - \frac{1}{\#(jk)} \Gamma^i_{jk} \vf{u^\alpha_{(jk)}} .
\]

A direct calculation confirms that these coordinate formul\ae\ for $p_\n$ are invariant under fibred coordinate transformations on $\pi$ and their prolongations to $\Jpia$ and $\Jbpi$; we also see that the restriction of $p_\n$ to the semiholonomic submanifold $\Jbpihat \subset \Jpia$ is the symmetrization map $\prm_2$, as we would expect. Note that $p_\n$ is a bundle morphism over $j^1\piao \to \piba$, not over $(\pia)_{1,0} \to \piba$:
\[
\begin{tikzcd}
\Jpia \arrow[swap]{d}{j^1\piao} \arrow[near start]{drr}{\piaao} && \Jbpi \arrow{d}{\pi_{2,1}} & &
T_{\Jbpi}\Jpia \arrow{r}{p_\n} \arrow[swap]{d}{Tj^1\piao} & T\Jbpi \arrow{d}{T\pi_{2,1}} \\
\Jpi \arrow[leftrightarrow]{rr} && \Jpi & & T\Jpi \arrow[leftrightarrow]{r} & T\Jpi
\end{tikzcd}
\]
We apply the same approach in the general case to obtain the infinitesimal projection
\[
p_\n \colon \ T_{\Jkpi}\Jpika \to T\Jkpi
\] with coordinate expressions
\begin{align*}
&p_\n \biggl( \vf{u^\alpha_{Ij}} \biggr) = \sum_{\abs{K}>0}^{k-\abs{I}} \frac{(I+K)!}{I! K!}
\Gamma^j_K \vf{u^\alpha_{I+K}}
\end{align*}
obtained by composing $\Asf$ with $S_\n$, with the necessary consequence that
\begin{align*}
&p_\n \biggl( \vf{u^\alpha_{J\cdot}} \biggr) = (1 - \abs{J}) \vf{u^\alpha_J}
- \sum_{I+1_j = J} \biggl( \sum_{\abs{K}=2}^{k-\abs{I}} \frac{(I+K)!}{I! K!} \Gamma^j_K \vf{u^\alpha_{I+K}} \biggr), \\
&p_\n \biggl( \vf{x^i} \biggr) = \vf{x^i} .
\end{align*}
This gives us the following result.

\begin{Theorem}
Let $\pi \colon E \to M$ be a fibred manifold, and let $\n$ be a symmetric linear connection on $M$. For each nonholonomic jet manifold $\Jpika$ there is a unique infinitesimal projection $p_\n$ satisfying $p_\n \circ T\iaka = \id_{T\Jkpi}$ and $p_\n|_{\Jkpihat} = \prm_k$, constructed by composing the isomorphism
\[
\Asf \colon \ V\pikaao \to \pikaa^* T^* M \otimes \pikaao^* V\pika
\]
arising from the affine structure of $\pikaao \colon \Jpika \to \Jkapi$ \big(restricted to points of $\Jkpi$\big) with the vertical tensor $S_\n$ on $\Jkpi$.
\end{Theorem}

\section{Homotopy operators for the horizontal differential}

We have seen that homotopy operators for the horizontal differential play an important part in the construction of Lepage equivalents satisfying the closure condition (and, indeed, of at most $1$-contact Lepage equivalents in general), and that locally such homotopy operators can be constructed using vertical endomorphisms. We have also noted that global homotopy operators (depending on a choice of a symmetric linear connection) have been shown to exist, but the construction in~\cite{And89} uses a quite different method, relating the horizontal differential on forms to an operator acting on evolutionary vector fields. It is therefore of some interest to see whether a global homotopy operator can be constructed directly for differential forms by using vertical endomorphisms. I conjecture that this can be done, and offer a possible method of doing so. The proposed formula has been checked for small values of the parameters $p$, $q$ and $r$ (see Appendix~\ref{appendixA} for an example calculation); although the general result might be amenable to a~direct calculation, there may well be a more geometric method of approaching it.

There are three ingredients in the proposed formula, which mimics the local formula described above. The vertical tensor $S_\n$, regarded as a map $\Omega^{p,q} \to \Xfk(M) \otimes \Omega^{p,q}$, has already been specified, and this can be iterated to give a map $S_\n^r \colon \Omega^{p,q} \to \odot^r \Xfk(M) \otimes \Omega^{p,q}$, where $\odot^r \Xfk(M)$ denotes the symmetric multivector fields on $M$. We shall also need a covariant version of the horizontal differential, which we shall denote $\dhn$; this will be a map $\odot^r \Xfk(M) \otimes \Omega^{p,q} \to \odot^r \Xfk(M) \otimes \Omega^{p,q+1}$, given on basis tensors by
\[
\dhn(X \otimes \omega) = \n X \wedge \omega + X \otimes \dh \omega
\]
and extended by multilinearity, symmetry and skewsymmetry. The final ingredient will be an operator $\Csf \colon \odot^r \Xfk(M) \otimes \Omega^{p,q} \to \odot^{r-1} \Xfk(M) \otimes \Omega^{p,q-1}$ contracting a vector component with a form component, again taking advantage of symmetry and skewsymmetry. The proposed homotopy operator is then $P_\n \colon \Omega^{p,q} \to \Omega^{p,q-1}$, where
\[
P_\n \omega = \sum_{r=0}^\infty \frac{(-1)^r (m-q)!}{p(m-q+r+1) r!}
\bigl( \Csf \circ \dhn \bigr)^r \Csf \bigl( S_\n^{r+1} \omega \bigr) .
\]

\section{Discussion}

One of the features of the approach taken in this paper is that it combines the use of finite and infinite jets. Variational problems are by their nature of finite order, and the various differential forms involved in their analysis are normally defined on a finite order jet manifold. Indeed, as we have seen, the properties of Lepage equivalents of first order and second order Lagrangians are rather different from those of higher order Lagrangians.

On the other hand, the variational bicomplex is best considered on the infinite jet manifold. In~\cite{And89}, a subcomplex called the Jacobian subcomplex which is projectable to a finite order jet manifold is shown after lengthy calculations to be locally exact; but no mention is made of a~homotopy operator acting on forms which are not $\dh$-closed. The operators $\Phat$ and $\Ptil$ described earlier, although acting on all the forms on each finite order jet manifold, generally increase their order. It seems to be the case that the complexity of ascertaining a bound on the order of the forms obscures the homotopy structure of the problem, and indeed the potential for a global solution. The alternative approach in~\cite{Voi22}, which involves a single homotopy operator for the variational derivative (and thus, essentially, for the vertical differential) avoids this problem, but then cannot reduce to the classical fundamental Lepage equivalent for first order Lagrangians; in addition, global versions are likely to be constrained by topological considerations.

The investigations in the second half of the paper suggest that vertical endomorphisms, when glued together as a vertical tensor using a symmetric linear connection, could be a significant part of the geometry of the jet bundle structure on a fibred manifold. If the conjecture that they define a global homotopy operator for $\dh$ is correct, then the simple formula
\[
\vartheta_{\lambda, \n} = \sum_{p=0}^m (-P_\n \dv)^p \lambda
\]
will give a Lepage equivalent of the Lagrangian $\lambda$ satisfying the closure property without the need for a separate choice of $\vartheta_\lambda^{(1)}$ to start the recursion. There will, though, be the question of whether the truncated form $\lambda - (P_\n \dv) \lambda$ is the same as the Poincar\'{e}--Cartan form constructed using the infinitesimal projections $p_\n$.

A final observation is that we have not explicitly addressed the question of whether it is possible to find, for second order Lagrangians, a geometrical construction of a Lepage equivalent satisfying the closure condition independently of any connection, as can be done for the Poincar\'{e}--Cartan form and the Carath\'{e}odory form. I suspect that this will not be the case.

\appendix

\section{An example calculation}\label{appendixA}

We consider the form $\omega = f^i_{\alpha m} {\rm d}x^m \otimes \theta^\alpha_i$, where $p = q = 1$ and the form is projectable to $\Jbpi$, so that the formula is
\[
P_\n \omega = \sum_{r=0}^\infty \frac{(-1)^r (m-1)!}{(m+r) r!} \bigl( \Csf \circ \dhn \bigr)^r \Csf \bigl( S_\n^{r+1} \omega \bigr) .
\]
We obtain
\[
\dh \omega = \bigl(d_l f^i_{\alpha m}\bigr) {\rm d}x^l \wedge {\rm d}x^m \wedge \theta^\alpha_i + f^i_{\alpha m} {\rm d}x^l \wedge {\rm d}x^m \wedge \theta^\alpha_{(il)}
\]
so that
\begin{gather*}
S_\n (\dh \omega) = \bigl(d_l f^i_{\alpha m}\bigr) \p_i \otimes {\rm d}x^l \wedge {\rm d}x^m \wedge \theta^\alpha
+ f^i_{\alpha m} \Gamma^k_{il} \p_k \otimes {\rm d}x^l \wedge {\rm d}x^m \wedge \theta^\alpha
\\ \hphantom{S_\n (\dh \omega) =}{}
+ f^h_{\alpha m} \p_k \otimes {\rm d}x^k \wedge {\rm d}x^m \wedge \theta^\alpha_h
+ f^k_{\alpha m} \p_k \otimes {\rm d}x^h \wedge {\rm d}x^m \wedge \theta^\alpha_h , \\
\Csf S_\n (\dh \omega) = \bigl(d_i f^i_{\alpha j}\bigr) {\rm d}x^j \wedge \theta^\alpha - \bigl(d_j f^i_{\alpha i}\bigr) {\rm d}x^j \wedge \theta^\alpha
+ f^i_{\alpha j} \Gamma^k_{ik} {\rm d}x^j \wedge \theta^\alpha - f^i_{\alpha k} \Gamma^k_{ij} {\rm d}x^j \wedge \theta^\alpha
\\ \hphantom{\Csf S_\n (\dh \omega) =}{}
 + m \omega - f^i_{\alpha i} {\rm d}x^j \wedge \theta^\alpha_j
\end{gather*}
and
\begin{gather*}
S_\n^2 (\dh \omega) = 2 \p_i \otimes \p_l \otimes \bigl( f^i_{\alpha m} {\rm d}x^l \wedge {\rm d}x^m \wedge \theta^\alpha \bigr) , \\
\Csf S_\n^2 (\dh \omega) = 2 m \p_i \otimes \bigl( f^i_{\alpha j} {\rm d}x^j \wedge \theta^\alpha \bigr)
- 2 \p_j \otimes \bigl( f^i_{\alpha i} {\rm d}x^j \wedge \theta^\alpha \bigr) , \\
\tfrac{1}{2} \dhn \Csf S_\n^2 \dh \omega
 = m \Gamma^k_{ih} \p_k \otimes {\rm d}x^h \wedge \bigl( f^i_{\alpha j} {\rm d}x^j \wedge \theta^\alpha \bigr)
+ m \p_i \otimes \bigl( \bigl(d_k f^i_{\alpha j}\bigr) {\rm d}x^k \wedge {\rm d}x^j \wedge \theta^\alpha \bigr)
\\ \hphantom{\tfrac{1}{2} \dhn \Csf S_\n^2 \dh \omega=}{}
 + m \p_i \otimes \bigl( f^i_{\alpha j} {\rm d}x^k \wedge {\rm d}x^j \wedge \theta^\alpha_k \bigr)
- \Gamma^k_{jh} \p_k \otimes {\rm d}x^h \wedge \bigl( f^i_{\alpha i} {\rm d}x^j \wedge \theta^\alpha \bigr)
\\ \hphantom{\tfrac{1}{2} \dhn \Csf S_\n^2 \dh \omega=}{}
 - \p_j \otimes \bigl( \bigl(d_k f^i_{\alpha i}\bigr) {\rm d}x^k \wedge {\rm d}x^j \wedge \theta^\alpha \bigr)
- \p_j \otimes \bigl( f^i_{\alpha i} {\rm d}x^k \wedge {\rm d}x^j \wedge \theta^\alpha_k \bigr) , \\
\tfrac{1}{2} \Csf \dhn \Csf S_\n^2 (\dh \omega)
 = m \Gamma^k_{ik} f^i_{\alpha j} {\rm d}x^j \wedge \theta^\alpha
- m \Gamma^k_{ij} {\rm d}x^j \wedge f^i_{\alpha k} \theta^\alpha
\\ \hphantom{\tfrac{1}{2} \Csf \dhn \Csf S_\n^2 (\dh \omega)=}{}
 + m \bigl(d_i f^i_{\alpha j}\bigr) {\rm d}x^j \wedge \theta^\alpha + m \omega
- \bigl(d_j f^i_{\alpha i}\bigr) {\rm d}x^j \wedge \theta^\alpha - f^i_{\alpha i} {\rm d}x^j \wedge \theta^\alpha_j .
\end{gather*}
On the other hand,
\begin{align*}
&S_\n \omega = \p_i \otimes \bigl(f^i_{\alpha m} {\rm d}x^m \wedge \theta^\alpha \bigr) , \\
&\Csf S_\n \omega = f^i_{\alpha i} \theta^\alpha , \\
&\dh (\Csf S_\n \omega) = \bigl(d_j f^i_{\alpha i}\bigr) {\rm d}x^j \wedge \theta^\alpha + f^i_{\alpha i} {\rm d}x^j \wedge \theta^\alpha_j,
\end{align*}
so that
\begin{align*}
\tfrac{1}{2} \Csf \dhn \Csf S_\n^2 (\dh \omega) + \dh (\Csf S_\n \omega)
={}& m \Gamma^k_{ik} f^i_{\alpha j} {\rm d}x^j \wedge \theta^\alpha - m \Gamma^k_{ij} {\rm d}x^j \wedge f^i_{\alpha k} \theta^\alpha \\
&{}+ m \bigl(d_i f^i_{\alpha j}\bigr) {\rm d}x^j \wedge \theta^\alpha + m \omega .
\end{align*}
But from
\begin{align*}
&\Csf S_\n (\dh \omega) = \bigl(d_i f^i_{\alpha j}\bigr) {\rm d}x^j \wedge \theta^\alpha
+ f^i_{\alpha j} \Gamma^k_{ik} {\rm d}x^j \wedge \theta^\alpha - f^i_{\alpha k} \Gamma^k_{ij} {\rm d}x^j \wedge \theta^\alpha + m \omega - \dh (\Csf S_\n \omega)
\end{align*}
we see that
\[
\tfrac{1}{2} \Csf \dhn \Csf S_\n^2 (\dh \omega) = m \Csf S_\n \dh \omega - m(m-1) \omega + (m-1) \dh \Csf S_\n \omega,
\]
so that
\[
\omega = \biggl( \frac{1}{m-1} \Csf S_\n - \frac{1}{2m(m-1)} \Csf \dhn \Csf S_\n^2 \biggr) \dh \omega
+ \dh \biggl( \frac{1}{m} \Csf S_n \omega \biggr) .
\]

\subsection*{Acknowledgements}

I should like to acknowledge correspondence with Nicoleta~Voicu which encouraged me to return to this topic after a number of years. Some results from this paper were presented at a meeting in Torino in honour of Marco Ferraris in June 2023, and at the International Summer School on Global Analysis and Applications in Pre\v{s}ov in August 2023.

\pdfbookmark[1]{References}{ref}
\LastPageEnding

\end{document}